\documentclass{article}%
\usepackage{amsmath}
\usepackage{amsfonts}
\usepackage{amssymb}
\usepackage{graphicx}%
\newtheorem{theorem}{Theorem}

\newtheorem{definition}[theorem]{Definition}

\newtheorem{remark}[theorem]{Remark}

\begin{document}

\title{{\LARGE Spectral problems in inhomogeneous media, spectral parameter power
series and transmutation operators}}
\author{Vladislav V. Kravchenko and Sergii M. Torba\\{\small Department of Mathematics, CINVESTAV del IPN, Unidad Queretaro, }\\{\small Libramiento Norponiente No. 2000, Fracc. Real de Juriquilla,
Queretaro, }\\{\small Qro. C.P. 76230 MEXICO, e-mail: vkravchenko@qro.cinvestav.mx}}
\maketitle

\begin{abstract}
We give a brief overview of recent developments in Sturm-Liouville theory
concerning operators of transmutation (transformation) and spectral parameter
power series (SPPS) and propose a new method for numerical solution of
corresponding spectral problems.

\end{abstract}

\section{Spectral parameter power series (SPPS)}

Let $f\in C^{2}(a,b)\cap C^{1}[a,b]$ be a complex valued function and
$f(x)\neq0$ for any $x\in\lbrack a,b]$. The interval $(a,b)$ is assumed being
finite. Let $x_{0}\in\lbrack a,b]$. We introduce the infinite system of
functions $\left\{  \varphi_{k}\right\}  _{k=0}^{\infty}$ defined as follows
\begin{equation}
\varphi_{k}(x)=%
\begin{cases}
f(x)X^{(k)}(x), & k\text{\ odd,}\\
f(x)\widetilde{X}^{(k)}(x), & k\text{\ even,}%
\end{cases}
\label{phik}%
\end{equation}
where
\[
X^{(0)}(x)\equiv1,\qquad X^{(n)}(x)=n\int_{x_{0}}^{x}X^{(n-1)}(s)\left(
f^{2}(s)\right)  ^{(-1)^{n}}\,\mathrm{d}s,\quad n=1,2,\ldots
\]
and
\[
\widetilde{X}^{(0)}\equiv1,\qquad\widetilde{X}^{(n)}(x)=n\int_{x_{0}}%
^{x}\widetilde{X}^{(n-1)}(s)\left(  f^{2}(s)\right)  ^{(-1)^{n-1}}%
\,\mathrm{d}s,\quad n=1,2,\ldots.
\]
Together with the system of functions (\ref{phik}) we define the functions
$\{\psi_{k}\}_{k=0}^{\infty}$ using the \textquotedblleft second
half\textquotedblright\ of the recursive integrals,
\begin{equation*}
\psi_{k}(x)=%
\begin{cases}
\dfrac{\widetilde{X}^{(k)}(x)}{f(x)}, & k\text{\ odd,}\\
\dfrac{X^{(k)}(x)}{f(x)}, & k\text{\ even.}%
\end{cases}
\end{equation*}

As we show below the introduced families of functions are closely related to
the one-dimensional Schr\"{o}dinger equation of the form $u^{\prime\prime
}-qu=\lambda u$ where $q$ is a complex-valued continuous function. Slightly
more general families of functions can be studied in relation to
Sturm-Liouville equations of the form $(py^{\prime})^{\prime}+qy=\lambda ry$.
Their definition based on a corresponding recursive integration procedure is
given in \cite{APFT}, \cite{KrPorter2010}, \cite{KKRosu}. The system
(\ref{phik}) is closely related to the notion of the $L$-basis introduced and
studied in \cite{Fage}.

The following result obtained in \cite{KrCV08} (for additional details and
simpler proof see \cite{APFT} and \cite{KrPorter2010}) establishes the
relation of the system of functions $\left\{  \varphi_{k}\right\}
_{k=0}^{\infty}$ and $\left\{  \psi_{k}\right\}  _{k=0}^{\infty}$ to the
Sturm-Liouville equation.

\begin{theorem}
\label{ThGenSolSturmLiouville} Let $q$ be a continuous complex valued function
of an independent real variable $x\in\lbrack a,b]$ and $\lambda$ be an
arbitrary complex number. Suppose there exists a solution $f$ of the equation
\begin{equation}
f^{\prime\prime}-qf=0 \label{SLhom}%
\end{equation}
on $(a,b)$ such that $f\in C^{2}(a,b)\cap C^{1}[a,b]$ and $f(x)\neq0$\ for any
$x\in\lbrack a,b]$. Then the general solution $u\in C^{2}(a,b)\cap C^{1}[a,b]$
of the equation
\begin{equation}
u^{\prime\prime}-qu=\lambda u \label{SLlambda}%
\end{equation}
on $(a,b)$ has the form $u=c_{1}u_{1}+c_{2}u_{2}$ where $c_{1}$ and $c_{2}$
are arbitrary constants,
\begin{equation}
u_{1}=\sum_{k=0}^{\infty}\frac{\lambda^{k}}{(2k)!}\varphi_{2k}\qquad
\text{and}\qquad u_{2}=\sum_{k=0}^{\infty}\frac{\lambda^{k}}{(2k+1)!}%
\varphi_{2k+1} \label{u1u2}%
\end{equation}
and both series converge uniformly on $[a,b]$ together with the series of the
first derivatives which have the form%
\begin{multline*}
u_{1}^{\prime}=f^{\prime}+\sum_{k=1}^{\infty}\frac{\lambda^{k}}{(2k)!}\left(
\frac{f^{\prime}}{f}\varphi_{2k}+2k\,\psi_{2k-1}\right)  \qquad\text{and}%
\\
u_{2}^{\prime}=\sum_{k=0}^{\infty}\frac{\lambda^{k}}{(2k+1)!}\left(
\frac{f^{\prime}}{f}\varphi_{2k+1}+\left(  2k+1\right)  \psi_{2k}\right)  .
\end{multline*}
The series of the second derivatives converge uniformly on any segment
$[a_{1},b_{1}]\subset(a,b)$.
\end{theorem}

The representation (\ref{u1u2}) offers the linearly independent solutions of
(\ref{SLlambda}) in the form of spectral parameter power series (SPPS). The
way of how the expansion coefficients in (\ref{u1u2}) are calculated via the
recursive integration is relatively simple and straightforward, this is why
the estimation of the rate of convergence of the series (\ref{u1u2}) presents
no difficulty, see \cite{KrPorter2010}. Moreover, in \cite{CamposKr} a
discrete analogue of Theorem \ref{ThGenSolSturmLiouville} was established and
the discrete analogues of the series (\ref{u1u2}) resulted to be finite sums.

\begin{remark}
\label{RemInitialValues}It is easy to see that the solutions $u_{1}$ and
$u_{2}$ defined by (\ref{u1u2}) satisfy the following initial conditions
\begin{align*}
u_{1}(x_{0}) &  =f(x_{0}), & u_{1}^{\prime}(x_{0}) &  =f^{\prime}(x_{0}),\\
u_{2}(x_{0}) &  =0, & u_{2}^{\prime}(x_{0}) &  =1/f(x_{0}).
\end{align*}

\end{remark}

\begin{remark}
\label{RemarkNonVanish} It is worth mentioning that in the regular case the
existence and construction of the required $f$ presents no difficulty. Let $q$
be real valued and continuous on $[a,b]$. Then (\ref{SLhom}) possesses two
linearly independent regular solutions\/ $v_{1}$ and $v_{2}$ whose zeros
alternate. Thus one may choose $f=v_{1}+iv_{2}$. Moreover, for the
construction of $v_{1}$ and $v_{2}$ in fact the same SPPS method may be used
\cite{KrPorter2010}.
\end{remark}

The SPPS representation (\ref{u1u2}) for solutions of the Sturm-Liouville
equation (\ref{SLlambda}) is very convenient for writing down the dispersion
(characteristic) relations in an analytical form. This fact was used in
\cite{BKK}, \cite{CKOR}, \cite{KKRosu}, \cite{KiraRosu2010},
\cite{KrPorter2010}, \cite{KV2011} for approximating solutions of different
spectral and scattering problems.

\section{Transmutation operators}

Let $E$ be a linear topological space and $E_{1}$ its linear subspace (not
necessarily closed). Let $A$ and $B$ be linear operators: $E_{1}\rightarrow E$.

\begin{definition}
A linear invertible operator $T$ defined on the whole $E$ such that $E_{1}$ is
invariant under the action of $T$ is called a transmutation operator for the
pair of operators $A$ and $B$ if it fulfills the following two conditions.

\begin{enumerate}
\item Both the operator $T$ and its inverse $T^{-1}$ are continuous in $E$;

\item The following operator equality is valid
\begin{equation}
AT=TB \label{ATTB}%
\end{equation}
or which is the same
\[
A=TBT^{-1}.
\]

\end{enumerate}
\end{definition}

Very often in literature the transmutation operators are called the
transformation operators. Here we keep ourselves to the original term coined
by Delsarte and Lions \cite{DelsarteLions1956}. Our main interest concerns the
situation when $A=-\frac{d^{2}}{dx^{2}}+q(x)$, $B=-\frac{d^{2}}{dx^{2}}$,
\ and $q$ is a continuous complex-valued function. Hence for our purposes it
will be sufficient to consider the functional space $E=C[a,b]$ with the
topology of uniform convergence and its subspace $E_{1}$ consisting of
functions from $C^{2}\left[  a,b\right]  $. One of the possibilities to
introduce a transmutation operator on $E$ (see, e.g., \cite{Marchenko})
consists in constructing a Volterra integral operator corresponding to a
midpoint of the segment of interest. Then it is convenient to consider a
symmetric segment $[-a,a]$ and hence the functional space $E=C[-a,a]$. It is
worth mentioning that other well known ways to construct the transmutation
operators (see, e.g., \cite{LevitanInverse}, \cite{Trimeche}) imply imposing
initial conditions on the functions and consequently lead to transmutation
operators satisfying (\ref{ATTB}) only on subclasses of $E_{1}$.

Thus, we consider the space $E=C[-a,a]$ and an operator of transmutation for
the defined above $A$ and $B$ can be realized in the form (see, e.g.,
\cite{LevitanInverse} and \cite{Marchenko}) of a Volterra integral operator
\begin{equation*}
Tu(x)=u(x)+\int_{-x}^{x}K(x,t)u(t)dt
\end{equation*}
where the kernel $K(x,t)$ is defined as a solution of a Goursat problem
\cite[Chapter 1]{Marchenko}. If the potential $q$ is $n$ times continuously
differentiable, the kernel $K(x,t)$ is $n+1$ times continuously differentiable
with respect to both independent variables.

In \cite{CKT} a parametrized family of operators $\mathbf{T}_{h}$,
$h\in\mathbb{C}$ was introduced, given by the integral expression
\begin{equation}
\mathbf{T}_{h}u(x)=u(x)+\int_{-x}^{x}\mathbf{K}(x,t;h)u(t)dt \label{Tmain}%
\end{equation}
where
\begin{equation}
\mathbf{K}(x,t;h)=\frac{h}{2}+K(x,t)+\frac{h}{2}\int_{t}^{x}%
\big(K(x,s)-K(x,-s)\big)\,ds. \label{Kmain}%
\end{equation}

The operator $\mathbf{T}_{h}$ maps a solution $v$ of an equation
$v^{\prime\prime}+\omega^{2}v=0$, where $\omega$ is a complex number, into a
solution $u$ of the equation $u^{\prime\prime}-q(x)u+\omega^{2}u=0$ with the
following correspondence of the initial values
\begin{equation}
u(0)=v(0),\qquad u^{\prime}(0)=v^{\prime}(0)+hv(0). \label{Th Initial Values}%
\end{equation}

\begin{theorem}
\cite{KTObzor} Suppose the potential $q$ satisfies either of the following two
conditions: 1) $q\in C^{1}[-a,a]$; 2) $q\in C[-a,a]$ and there exists a
particular complex-valued solution $g$ of (\ref{SLhom}) non-vanishing on
$[-a,a]$. Then the operator $\mathbf{T}_{h}$ given by \eqref{Tmain} satisfies
the equality
\begin{equation*}
\left(  -\frac{d^{2}}{dx^{2}}+q(x)\right)  \mathbf{T}_{h}[u]=\mathbf{T}%
_{h}\left[  -\frac{d^{2}}{dx^{2}}(u)\right]  
\end{equation*}
for any $u\in C^{2}[-a,a]$.
\end{theorem}

Suppose now that a function $f$ is a solution of \eqref{SLhom}, non-vanishing
on $[-a,a]$ and normalized as $f(0)=1$. Let $h:=f^{\prime}(0)$ be some complex
constant. Define as before the system of functions $\{\varphi_{k}%
\}_{k=0}^{\infty}$ by this function $f$ and by \eqref{phik}. The following
theorem states that the operator $\mathbf{T}_{h}$ transmutes powers of $x$
into the functions $\varphi_{k}$.

\begin{theorem}
\cite{CKT}\label{Th Transmutation of Powers} Let $q$ be a continuous complex
valued function of an independent real variable $x\in\lbrack-a,a]$, and $f$ be
a particular solution of \eqref{SLhom} such that $f\in C^{2}(-a,a)$ together
with $1/f$ are bounded on $[-a,a]$ and normalized as $f(0)=1$, and let
$h:=f^{\prime}(0)$, where $h$ is a complex number. Then the operator
\eqref{Tmain} with the kernel defined by \eqref{Kmain} (with this particular
$h$) transforms $x^{k}$ into $\varphi_{k}(x)$ for any $k\in\mathbb{N}_{0}$.
\end{theorem}

Thus, the system of functions $\{\varphi_{k}\}$ may be obtained as the result
of the Volterra integral operator acting on powers of the independent
variable. As was mentioned before, this offers an algorithm for transmuting
functions in situations when the explicit form of $\mathbf{K}(x,t;h)$ is
unknown. Moreover, properties of the Volterra integral operator such as
boundedness and bounded invertibility in many functional spaces gives us a
tool to prove the completeness of the system of function $\{\varphi_{k}\}$ in
various situations. For example, \textbf{the system }$\{\varphi_{k}%
\}_{k=0}^{\infty}$\textbf{ is complete in }$C[-a,a]$\textbf{.}

\section{A new method for solving spectral problems}

Theorem \ref{Th Transmutation of Powers} together with
(\ref{Th Initial Values}) suggests the following approach for solving spectral
problems for the Sturm-Liouville equation (\ref{SLlambda}). Here as an example
we consider the Sturm-Liouville problem with the following boundary conditions
though obviously the method can be applied in a much more general situation.
Thus, consider the equation
\begin{equation}
u^{\prime\prime}-qu=-\beta^{2}u\label{SLwithBeta}%
\end{equation}
($-\beta^{2}=\lambda$) with the conditions
\begin{equation}
u(0)=u(1)=0.\label{Dirichlet conditions}%
\end{equation}
It is well known that an almost optimal uniform approximation of the functions
sine and cosine on the segment $[-1,1]$ by polynomials is achieved using the
Chebyshev polynomials. The corresponding representations have the form
\cite[p. 104]{Suetin}%
\begin{align}
\sin\beta x & =2\sum_{m=0}^{\infty}(-1)^{m}J_{2m+1}(\beta)T_{2m+1}(x),\\
\cos\beta x & =J_{0}(\beta)+2\sum_{m=0}^{\infty}(-1)^{m}J_{2m}(\beta
)T_{2m}(x)\label{sincos}%
\end{align}
where $J_{n}$ denotes the Bessel function of the first kind and $T_{n}$ is the
Chebyshev polynomial of $n$-th order, $T_{n}(x)=\cos(n\arccos x)$.

From (\ref{Th Initial Values}) we have that $u(x)=\mathbf{T}_{h}\left(
\sin\beta x\right)  $ is a solution of (\ref{SLwithBeta}) satisfying the first
of the conditions (\ref{Dirichlet conditions}). Thus, $\beta^{2}$ is an
eigenvalue of (\ref{SLwithBeta}), (\ref{Dirichlet conditions}) iff
$\mathbf{T}_{h}\left(  \sin\beta x\right)  $ at the point $x=1$ equals zero.
Notice that
\begin{equation}
\begin{split}
\mathbf{T}_{h}\left(  \sin\beta x\right)  =&2\sum_{m=0}^{\infty}(-1)^{m}%
J_{2m+1}(\beta)\mathbf{T}_{h}\left(  T_{2m+1}(x)\right)  \\
\cong&2\sum_{m=0}%
^{N}(-1)^{m}J_{2m+1}(\beta)\mathbf{T}_{h}\left(  T_{2m+1}(x)\right).
\end{split}\label{ThSin}%
\end{equation}
In the last expression it is not difficult to calculate $\mathbf{T}_{h}\left(
T_{2m+1}(x)\right)  $ with a very good accuracy applying theorem
\ref{Th Transmutation of Powers}. Every $\mathbf{T}_{h}\left(  T_{2m+1}%
(x)\right)  $ is a linear combination of the functions $\varphi_{k}$,
$k=0,\ldots2m+1$ and the coefficient corresponding to $\varphi_{k}$ is the
same as the coefficient corresponding to $x^{k}$. Thus, computation of the
images of the Chebyshev polynomials under the action of the transmutation
operator $\mathbf{T}_{h}$ does not represent any difficulty. Taking into
account a very fast convergence of the series in (\ref{sincos}) one can obtain
a very good accuracy in the approximation of $\mathbf{T}_{h}\left(  \sin\beta
x\right)  $ with a relatively small $N$ in (\ref{ThSin}).

Thus, the proposed method for solving (\ref{SLwithBeta}),
(\ref{Dirichlet conditions}) consists of the following steps: 1) compute a
solution $f$ of (\ref{SLhom}) satisfying the conditions of Theorem
\ref{ThGenSolSturmLiouville}; 2) compute $N+1$ functions $\varphi_{k}$; 3)
taking the coefficients of $x^{k}$, $k=0,\ldots2m+1$ from the expression of
$T_{2m+1}$ compute $\mathbf{T}_{h}\left(  T_{2m+1}\right)  $ evaluated at
$x=1$; 4) find zeros of $\sum_{m=0}^{N}(-1)^{m}J_{2m+1}(\beta)\mathbf{T}%
_{h}\left(  T_{2m+1}\right)  \mid_{x=1}$.

\textbf{Example }$q(x)=e^{x}$, $u(0)=u(\pi)=0$. First, making an obvious
change of variables we reduce the problem to the interval $(0,1)$ and then
compare our results with the eigenvalues from \cite[Appendix A]{Pryce}. Taking
$N=18$ we obtain the results reported in the following table. The recursive integration involved in the construction of the functions $\varphi_k$ required in \eqref{ThSin} was performed by means of the Spline Toolbox of Matlab. On each step the integrand was converted into a spline and then integrated using the command \textsf{fnint}. A surprising accuracy with a very small number of functions participating in \eqref{ThSin} was achieved in all the considered numerical tests. Here due to the restrictions in space we illustrate this with the only example in which $N=18$. A phenomenon which we observed but up to now have not found an explanation is a surprisingly good approximation of higher order eigenvalues as can be seen from the table. Apparently the expression \eqref{ThSin} not only approximates well the solution of the equation for relatively small values of $\beta$ but also evaluated at $x=1$ contains information on the asymptotics of the eigenvalues of the problem. It is worth mentioning that on a usual netbook computer the computation of the eigenvalues (see the table below) implemented in the Matlab  takes only several seconds.

\noindent
\begin{center}\begin{tabular}
[c]{|c|c|c|}\hline
Order of the eigenvalue & Eigenvalues from \cite{Pryce} & Computed
eigenvalues\\\hline
1 & 4.8966693800 & 4.896669377\\\hline
2 & 10.045189893 & 10.045189883\\\hline
3 & 16.019267250 & 16.019267268\\\hline
4 & 23.266270940 & 23.266270969\\\hline
5 & 32.26370704 & 32.26370710\\\hline
6 & 43.2200196 & 43.2200184\\\hline
7 & 56.18159 & 56.18161\\\hline
8 & 71.15299 & 71.15255\\\hline
9 & 88.1321 & 88.1386\\\hline
10 & 107.11 & 107.05\\\hline
11 & 128.10 & 128.48\\\hline
17 & 296.07 & 296.52\\\hline
28 & 791.05 & 790.99\\\hline
43 & 1856.05 & 1856.53\\\hline
50 & 2507 & 2500\\\hline
\end{tabular}
\end{center}

\section*{Acknowledgements}
The work was partially supported by Conacyt via the project 166141.

\end{document}